\documentclass[12pt]{amsart}
\input{amssym.def}
\usepackage{graphicx}

\newtheorem{theorem}{Theorem}[section]                      
\newtheorem{proposition}{Proposition}[section]                       
\newtheorem{corollary}{Corollary}[section]  
\newtheorem{lemma}{Lemma}[section]              
\newtheorem{remark}{Remark}[section]
\newtheorem{question}{Question}
\newtheorem{definition}{Definition}[section]
\theoremstyle{remark} 
\newtheorem{example}{\bf Example}[section]

\newcommand{\s}{\vspace{0.3cm}}

\input epsf

\begin{document}

\title[Noded Teichm\"uller Spaces]
{Noded Teichm\"uller Spaces}

\author{Rub\'en A. Hidalgo and Alexander Vasil'ev}
\thanks{Partially supported by 
Projects Fondecyt 1030252, 1030373, 1040333 and Projects UTFSM 12.05.21, 
12.05.23}
\email{ruben.hidalgo@usm.cl}
\email{alexander.vasiliev@usm.cl and alexander.vasiliev@uib.no}
\address{Departamento de Matem\'atica, UTFSM, Casilla 110-V Valparaiso, Chile}
\address{Department of Mathematics, University of Bergen, Johannes Brunsgate 12, Bergen 5008, Norway}
\dedicatory{In memory of our colleague Miguel Bl\'azquez}

\begin{abstract}  
Let $G$ be a finitely generated Kleinian group and let $\Delta$ be
an invariant collection of components in its region of
discontinuity.  The Teichm\"uller space $T(\Delta,G)$ supported in $\Delta$, 
is the space of equivalence classes of quasiconformal
homeomorphisms with complex dilatation invariant under $G$ and supported in
$\Delta$. In this paper we propose a partial closure of
$T(\Delta,G)$ by considering certain deformations of the above homeomorphisms. Such a partial closure is denoted by $NT(\Delta,G)$ and called the {\it noded Teichm\"uller space} of $G$ supported in $\Delta$. Some concrete examples are discussed.
\end{abstract}

\maketitle

\section{Introduction}
An {\it analytically finite Riemann orbifold} $S$ is a closed Riemann surface of genus $\gamma$ together with a finite collection of conical points $x_{1},...,x_{n} \in S$ of orders $2 \leq v_{1},...,v_{n} \leq \infty$, respectively. It is said to have signature $(\gamma,n;v_{1},...,v_{n})$. In the case $n=3$ we say that $S$ is a {\it triangular orbifold}. If $n=0$, then $S$ is a closed Riemann surface, and if $n>0$ and $v_{j}=\infty$ far all $j$, then $S$ is an analytically finite punctured Riemann surface. Let $G$ be a Kleinian group with the region of discontinuity $\Omega$, and let $L \in G$ be a loxodromic transformation. We say that $L$ is {\it primitive} if it is not a nontrivial positive power of another loxodromic transformation in $G$. We say that $L$ is {\it simple loxodromic} if there is a simple arc on $\Omega$ which is invariant under $L$, we call it an {\it axis} of $L$, whose projection on $\Omega/G$ is a simple loop or a simple arc connecting two conical points of order $2$. A parabolic transformation $P \in G$ with a fixed point $p$ is called {\it double-cusped}, if any other parabolic element of $G$ commuting with $P$ belongs to the cyclic group generated by it and there are two tangent  open discs at $p$ in $\Omega$, which are invariant under the stabilizer of $p$ in $G$. Generalities on Kleinian groups may be found in Maskit's monograph \cite{M1}.

Let $G$ be a finitely generated Kleinian group $G$ with the region of discontinuity $\Omega(G)$ and limit set $\Lambda(G)=\widehat{\mathbb C}-\Omega(G)$. We proceed to recall the definition of the {\it Teichm\"uller space of $G$ supported at $\Delta$}. Details may be found in \cite{N}.
Let $\Delta \subset \Omega(G)$ be a collection of components of $\Omega(G)$ which is 
$G$-invariant, that is $\gamma (\Delta)=\Delta$ for every $\gamma \in G$. As a consequence of the Ahlfors finiteness theorem \cite{Ahlfors:finitud}, we have that $\Delta/G$ is a finite union of analytically finite Riemann orbifolds. Set 
$$L^{\infty}(\Delta,G)=
\{\mu \in L^{\infty}(\Delta); \; \mu(g(z)) \overline{g'(z)} = \mu(z) g'(z), \;
{\rm a.e.}\; \Delta,\; \mbox{for all} \; g \in G\},$$ and
$$L^{\infty}(\Delta,G)_{1}=\{\mu \in L^{\infty}(\Delta,G); \;
\|\mu\|_{\infty} < 1\}$$

The measurable functions $\mu$ in $L^{\infty}(\Delta,G)_{1}$ are called
{\it Beltrami coefficients} for $G$ supported in $\Delta$. 
An orientation preserving homeomorphism $w:\widehat{\mathbb C} \to
\widehat{\mathbb C}$ is called a {\it quasiconformal homeomorphism} if there is
a Beltrami coefficient $\mu \in L^{\infty}(\widehat{\mathbb C})_{1}$
(called the {\it complex dilatation} of $w$) such that $w$ has distributional partial
derivatives $\partial w$, $\overline{\partial} w$ in
$L^{2}_{loc}(\widehat{\mathbb C})$ satisfying the Beltrami equation
$$\overline{\partial}w(z)=\mu(z)\partial w(z), \quad \mbox{ a.e. $z \in
\widehat{\mathbb C}$ }$$

The following lemma is a 
classical result and its proof can be found, for instance, in \cite{L}.

\s
\noindent
\begin{lemma} Let 
$\mu \in L^{\infty}(\Delta,G)_{1}$ and let
$w:\widehat{\mathbb C} \to \widehat{\mathbb C}$
be a quasiconformal homeomorphism with complex dilatation $\mu$.
Then $w G w^{-1}$ is a Kleinian group with
region of discontinuity $w(\Omega(G))$.
\end{lemma}

\s

Let $G$ and $\Delta$ be defined as before and
let $\mu$ and $\nu$ be two Beltrami coefficients for $G$
supported in $\Delta$. If $w_{\mu}$ (respectively $w_{\nu}$) is a 
quasiconformal homeomorphism solving the Beltrami
equation for $\mu$ (respectively $\nu$), then we have a natural isomorphism
$\theta_{\mu}:G \to w_{\mu} G w_{\mu}^{-1}$ 
(respectively $\theta_{\nu}:G \to w_{\nu} G w_{\nu}^{-1}$).
We say that $\mu$ and $\nu$ are {\it Teichm\"uller equivalent}
if there is a M\"obius transformation $A$ satisfying the equality 
$\theta_{\mu}(g)=A \theta_{\nu}(g) A^{-1}$ for all $g \in G$.
In the case of a non-elementary $G$ equivalently,
$w_{\nu}$ and $A w_{\mu}$ coincide on the limit set of $G$. 
The {\it Teichm\"uller space} $T(\Delta,G)$ of
$G$  supported in $\Delta$ is the set of all Teichm\"uller equivalence classes of Beltrami coefficients for $G$
supported in $\Delta$. If $\Delta=\Omega(G)$, then we denote this space by
$T(G)$. In the particular case of a simply connected $\Delta$, the space $T(\Delta,G)$ is a copy of the Teichm\"uller space $T(S)$ of the Riemann orbifold $S=\Delta/G$, see \cite{N}.

If $G$ is non-elementary and
$\Delta/G$ is not the union of triangular orbifolds (we say that $G$ is {\it not
triangular}), then  $T(\Delta,G)$ is known to be a non-compact complex manifold of finite dimension (see, for instance, \cite{K-M} and \cite{N}). Let us write $\Delta/G=S_{1} \cup \cdots \cup S_{k}$, where $S_{j}$ ($j=1,...,k$) is an analytically finite Riemann orbifold.  Consider a maximal
collection $\Delta_{1}$,..., $\Delta_{k}$ of
non-equivalent components of $\Delta$, and set $G_{i}=\{g \in G :
g(\Delta_{i})=\Delta_{i}\}$. We may assume that  $\Delta_{i}/G_{i}=S_{i}$. As
a consequence of the results of \cite{Kra} and  \cite{S}, we have an
isomorphism $Q:T(\Delta,G) \to T(\Delta_{1},G_{1}) \times\cdots \times
T(\Delta_{k},G_{k})$. Moreover, the  results of Bers \cite{B3} and Maskit
\cite{M3}, assert that the universal cover of $T(\Delta,G)$ is the product of
the Teichm\"uller spaces $T(S_{1})\times \cdots \times T(S_{k})$. 
Let us denote by ${\mathcal M}_{i}$ the
moduli space of $S_{i}$, i.e., the space of conformal classes
of analytically finite Riemann orbifold structures on $S_{i}$, and let $\mathcal M$ stand for the disjoint union of ${\mathcal M}_{1}$,..., ${\mathcal
M}_{k}$. The holomorphic map $Q$ induces a natural holomorphic (branched) covering map  $\pi: T(\Delta,G) \to {\mathcal M}$.

It is known that ${\mathcal M}_{i}$ is compact only if $S_{i}$ is a Riemann orbifold with at most $3$ conical points. If we consider on $S_{i}$ a  finite collection of pairwise disjoint families of (i) simple closed geodesics disjoint from the conical points, and (ii) simple geodesic arcs connecting two conical points of order $2$, each of which does not contain other conical points, then, by a process of pinching them, we obtain a topological space called a {\it stable Riemann orbifold of type $S_{i}$}. We may talk on conformal homeomorphisms of stable Riemann orbifolds (orientation preserving homeomorphisms fixing the conical points with their orders, being holomorphic at the complement of the pinched portions, called {\it nodes}). In particular, we may talk on conformal classes of stable Riemann orbifolds of type $S_{i}$. A compactification $\widetilde{\mathcal M}_{i}$ of  ${\mathcal M}_{i}$, called the Deligne-Mumford's compactification, is obtained by adding  the conformal classes of stable Riemann orbifolds  of type $S_{i}$ to ${\mathcal M}_{i}$.  In this way, a compactification of $\mathcal M$ is obtained as the product of the above compactifications. 

At this point, it is natural to ask for a partial closure of  $T(\Delta,G)$, such that the border points of this closure correspond to some of the border points in Deligne-Mumford's compactification (as previously seen) in a natural sense. We shall construct such a clusure $NT(\Delta,G)$ in Section \ref{construccion}, so that the extra points that we add to $T(\Delta,G)$ produce, in terms of Kleinian groups, stable Riemann orbifolds in ${\mathcal M}$. The points in $NT(\Delta,G)-T(\Delta,G)$ correspond to deformations of the group $G$ by the process of approximation of double-cusped parabolic elements of $G$ by certain {\it primitive simple loxodromic} ones  (see the works of Keen, Series and Maskit in \cite{KMS}, and of Maskit in \cite{M2,M4}). The deformation is produced by some boundary points $\mu \in \partial L^{\infty}(G,\Delta)_{1}$, called {\it noded Beltrami differentials for $G$}. At the level of the Riemann orbifold $S=\Delta/G$ this means that we permit certain pairwise disjoint simple closed geodesics on $S$ and maybe some simple geodesic arcs connecting conical values of order $2$ to degenerate to points in order to produce a finite collection of noded Riemann orbifolds.  The loops and arcs which we consider in the degeneration process are the projections of appropriately chosen axes of the primitive simple loxodromic elements of the group that approach  double-cusped parabolic transformations. 

In the particular case of a Fuchsian group $G$ of the first kind and of
$\Delta$ as one of two discs in $\Omega(G)$, the constructed partial closure 
is  exactly Abikoff's Augmented Teichm\"uller space \cite{Ab1} \cite{Ab2}. If $G$ is a Schottky group and $\Delta=\Omega(G)$, then $NT(G,\Delta)$ coincides with the noded Schottky space as defined in \cite{H1}. In the survey \cite{G-W} is discussed some others closures of Schottky space. But, as long a group is geometrically finite and free of rank equal to the rank of $G$, then it belongs to the noded Schottky space $NT(G,\Delta)$ \cite{H1,H-M}. In Section \ref{construccion} we provide more examples to describe this partial closure and relate them to known partial closures.

If $G$ is a finitely generated Fuchsian group acting on
the upper-half plane ${\mathbb H}$, so that ${\mathbb H}/G$ is a closed Riemann
surface, then  a natural compactification of $T({\mathbb H},G)$, at the level
of quasiconformal deformation, due to Bers \cite{B2} is obtained by
the embedding of $T({\mathbb H},G)$  as a bounded domain into the space of quadratic 
holomorphic forms. Another natural compactification of the
Teichm\"uller space of a surface of genus $g$,  due
to Thurston, is obtained by adding measured projective laminations.
It is to large in our setup. 
Another compactification of the Teichm\"uller space can be found in \cite{M-P}.

\section{Noded Quasiconformal Deformations}
In this section, $G$ will denote a finitely generated Kleinian group and $\Delta \neq \emptyset$ a $G$-invariant collection of components of $\Omega(G)$. 

\subsection{Noded quasiconformal maps} We denote by $\overline{L^{\infty}(\Delta,G)_{1}}$ the closure of
$L^{\infty}(\Delta,G)_{1}$ inside $L^{\infty}(\Delta,G)$. For each $\mu$ from $\overline{L^{\infty}(\Delta,G)_{1}}$
we define the {\it region of discontinuity} of 
$\mu$, denoted by
$\Omega(\mu)$, as the set of all points $p \in \widehat{\mathbb C}$, 
such that there is 
an open neighborhood $U$ of $p \in U$ and $\|\mu|_{U}\|_{\infty} <1$.
Its complement $\Lambda(\mu)=\widehat{\mathbb C}-\Omega(\mu)$ is the {\it limit set}
of $\mu$. By  definition, the set $\Omega(\mu)$ is open and $\Lambda(\mu)$ is
compact. 

\s
\noindent
\begin{proposition} \label{Prop1}
For each $\mu \in 
\overline{L^{\infty}(\Delta,G)_{1}}$,  both $\Omega(\mu)$ and $\Lambda(\mu)$ are $G$-in\-va\-riant.
\end{proposition}
\begin{proof} It is enough to check the invariance of $\Omega(\mu)$.
Let $p \in \Omega(\mu)$, $U$ be an open set containing $p$, such that
$\|\mu|_{U}\|_{\infty} < 1$ and $g \in G$.  
Set $g(U)=V$. The invariance property $\mu(g(z))\overline{g'(z)}=\mu(z)g'(z)$
asserts that $\|\mu|_{V}\|_{\infty}<1$. 
In particular, $g(p) \in \Omega(\mu)$. Since $g^{-1}$ also belongs to $G$,
we have the invariance property of $\Omega(\mu)$ as claimed. 
\hspace{\fill}\end{proof}

\s

If $U \subset \widehat{\mathbb C}$ is an open set, then we denote by $L^{2,1}_{loc}(U)$ the complex vector space of maps $w:U \to \widehat{\mathbb C}$ with locally integrable distributional derivatives. 

\s
\noindent
\begin{definition}
Let $\mu \in \overline{L^{\infty}(\Delta,G)_{1}}$ and let $w \in L^{2,1}_{loc}(\Omega(\mu))$ be an orientation-preserving map. We say that $w$ is a {\it noded quasiconformal map} with dilatation $\mu$ if \begin{itemize}
\item[(i)] there is a component of $\Omega(\mu)$ homeomorphically mapped by $w$ onto its image;
\item[(ii)] $ \overline{\partial} w(z) = \mu(z) \partial w(z), \,\,\, \mbox{\rm a.e. $z \in 
\Omega(\mu)$}$;
\item[(iii)] there is a sequence $\mu_{n} \in L^{\infty}(\Delta,G)_{1}$,
converging to $\mu$ almost everywhere in $\Omega(\mu)$;
\item[(iv)] there is a sequence $w_{n}:\widehat{\mathbb C} \to \widehat{\mathbb C}$ of
quasiconformal homeomorphisms with complex dilatations $\mu_{n}$, converging to $w$
locally uniformly in $\Omega(\mu)$.

\end{itemize}
\end{definition}

\s

The following existence result is classical and we give the proof as a matter
of completeness. 

\s
\noindent
\begin{proposition}\label{Prop2}  
If $\mu \in \overline{L^{\infty}(\Delta,G)_{1}}$
is such that $\Omega(\mu) \neq \emptyset$, $\Omega_{1}$ is a connected component
of $\Omega(\mu)$, and $x_{1}$, $x_{2}$, $x_{3}$ are three different points
in $\Omega_{1}$, then there is a noded quasiconformal map $w$ with dilatation $\mu$  
fixing the points $x_{1}$, $x_{2}$ and $x_{3}$,
which is a homeomorphism when restricted to $\Omega_{1}$.
\end{proposition}
\begin{proof} Let us assume $\mu \in \overline{L^{\infty}(\Delta,G)_{1}} - L^{\infty}(\Delta,G)_{1}$.
By Proposition \ref{Prop1}, both sets $\Omega(\mu)$ and $\Lambda(\mu)$ are $G-$invariant. Let us consider the maximal collection $\Omega_{1}$, $\Omega_{2},\dots$ of connected components of $\Omega(\mu)$ which  are not $G-$equivalent.

We set domains $A_{n}^{i} \subset \Omega_{i}$, such that
\begin{itemize}
\item[(i)] $\overline{A_{n}^{i}} \subset A_{n+1}^{i}$;
\item[(ii)] $\cup_{n}A_{n}^{i} = \Omega_{i}$;
\item[(iii)] $\mu$ restricted to $\cup_{i}A_{n}^{i}$ has essential sup norm
less than one;
\item[(iv)] the points $x_{1}$, $x_{2}$, $x_{3}$ belong to $A_{1}^{1}$.
\end{itemize}

Set 
$$
\mu_{n}=\left\{\begin{array}{ll}
\mu, & \mbox{in} \;G(\cup_{i}A_{n}^{i})\\
0, & \mbox{otherwise}.\end{array}
\right.
$$

We have a sequence $\mu_{n} \in L^{\infty}(\Delta,G)_{1}$ and let  
$w_{n}:\widehat{\mathbb C} \to \widehat{\mathbb C}$ be  quasiconformal 
homeomorphisms with complex dilatations $\mu_{n}$, fixing $x_{1}$, $x_{2}$, $x_{3}$.
Let us fix $m$. Then, for $n \geq m$ the sequence of 
quasiconformal homeomorphisms $w_{n}$ restricted to $A_{m}^{1}$
are  $K_{m}-$quasiconformal homeomorphisms with
$K_{m}=\displaystyle{
\frac{\mbox{$1+\|\mu_{m}\|_{\infty}$}}{\mbox{$1-\|\mu_{m}\|_{\infty}$}}}$, 
each of which  fixes the points $x_{1}$, $x_{2}$, $x_{3}$.
 Results on quasiconformal maps (see  \cite[page 14]{L}) assert that there
is a subsequence of $\{w_{n}\}$  converging locally uniformly to a map
$f_{m}^{1}:A_{m}^{1} \to \widehat{\mathbb C}$.  Of course, $f_{m}^{1}$
fixes the points $x_{1}$, $x_{2}$, $x_{3}$. 
It follows from \cite[page 15]{L}, that
$f_{m}^{1}$ is a $K_{m}-$quasiconformal homeomorphism (onto its image, so $f_{m}^{1}$ is
an  embedding).
Similarly, we consider a subsequence of maps on the domain $A_{m}^{2}$, which is
a normal sequence in such a domain. One may find a
subsequence of the above one that converges locally uniformly in $A_{m}^{2}$
either to a constant map or to a quasiconformal homeomorphism $f_{m}^{2}$. 
We proceed 
inductively with the domains $A_{m}^{3}$, $A_{m}^{4}$,..., to get subsequences
converging locally uniformly to some constant or a quasiconformal 
homeomorphism $f_{m}^{i}$.
Applying the diagonal process, we get a subsequence of the original sequence
$w_{n}$ that converges locally uniformly in each domain $A_{m}^{i}$ to a 
map $f_{m}^{i}$.
We use this sequence to work similarly in each domain $A_{m+1}^{i}$.
In this way we get sequences with the required convergence property for each $m$.
We use again the diagonal process to obtain the desired 
sequence of maps. The map $w$ is given locally by the functions
$f_{m}^{i}$.
\hspace{\fill}
\end{proof}

\subsection{Noded Beltrami coefficients for $G$ supported in $\Delta$}
By the stereographic projection, we may see the Riemann sphere as the unit sphere in ${\mathbb R}^{3}$. This allows us to consider the spherical metric and work with the spherical diameter of a subset of $\widehat{\mathbb C}$. 

\s
\noindent
\begin{definition} \label{nodedfamily}
Let $\mu \in \overline{L^{\infty}(G,\Delta)_{1}}$.
A countable collection ${\mathcal F}_{\mu}=\{ \alpha_{1},\alpha_{2},...\}$ of pairwise disjoint simple arcs
(including end points)
is called a {\it noded family of arcs associated with $\mu$}, if the
following properties hold  
\begin{itemize}
\item[(1)] $\alpha_{n}^{*} \subset \Delta$, where $\alpha_{n}^{*}$ denotes
$\alpha_{n}$ minus both extremes; 
\item[(2)] the spherical diameter of
$\alpha_{n}$ goes to $0$ as $n$ goes to $\infty$; 
\item[(3)] $\Lambda(\mu)=\overline{\cup_{j=1}^{\infty}\alpha_{j}}$;
\item[(4)] $\Omega(\mu) \subset \widehat{\mathbb C}$ is a dense subset; 
\item[(5)] the group $G_{n}=\{g \in G; g(\alpha_{n})=\alpha_{n}\}$, is either
a cyclic loxodromic group or a ${\mathbb Z}_{2}-$extension of a cyclic 
loxodromic group.
\end{itemize}

\end{definition}

\s

We remark that there are Beltrami coefficients $\mu \in \overline{L^{\infty}(G,\Delta)_{1}}$ for which there is no associated noded family of arcs, for instance, if  $\mu \in L^{\infty}(G,\Delta)_{1}$. In Section \ref{ejemplo} we construct an example of $\mu \in \overline{L^{\infty}(G,\Delta)_{1}}$ with a noded family of arcs, see also Example \ref{Ex1}.

\s
\noindent
\begin{definition}\label{nodedbeltrami}
We define the set $L^{\infty}_{{\rm noded}}(\Delta,G)$
of {\it noded Beltrami coefficients for $G$ supported in $\Delta$} as
a set of those
$\mu \in  \overline{L^{\infty}(\Delta,G)_{1}}$
for which there is an associated noded family of arcs ${\mathcal F}_{\mu}=\{ \alpha_{1},\alpha_{2},...\}$, such that 
there is a continuous map $w:\widehat{\mathbb C} \to \widehat{\mathbb C}$, 
called a noded quasiconformal deformation of $G$, with the complex dilatation
$\mu$,  and that $w \in L^{2,1}_{loc}(\Omega(\mu))$, 
satisfying the following properties:  
\begin{itemize}
\item[(1)] $w$ is injective in $\widehat{\mathbb C}-{\mathcal F}_{\mu}$;  
\item[(2)] $w$, when restricted to $\Omega(\mu)$, is a noded quasiconfromal map with complex dilataton $\mu$; 
\item[(3)] the restriction of $w$ to each arc $\alpha_{i}$ is a constant 
$p_{i}$, where $p_{i} \neq p_{j}$ for $i \neq j$. 
\end{itemize}
\end{definition}

\s
\noindent
\begin{remark}
 Let $\mu \in L^{\infty}_{{\rm noded}}(\Delta,G)$ and let
$w:\widehat{\mathbb C} \to \widehat{\mathbb C}$ be a 
noded quasiconformal deformation of $G$ with the complex dilatation $\mu$. Then,
the following statements are true
\begin{itemize}
\item[(1)] if $\alpha_{n} \in {\mathcal F}_{\mu}$, then 
both end points are the fixed points
of a loxodromic element of $G$. Such a loxodromic element keeps 
the connected component of $\Omega(G)$ containing $\alpha_{n}$
invariant;
\item[(2)] if $\Lambda(\mu) \neq \emptyset$, then 
$\Lambda(G) \subset \Lambda(\mu)$. This is a consequence of Proposition E.4 in \cite[page 96]{M1};
\item[(3)] $w(\Omega(u)) \cap w(\Lambda(\mu))=\emptyset$. Indeed, if there
were points $p_{1} \in \Omega(\mu)$ and $p_{2} \in \Lambda(\mu)$, such that
$w(p_{1})=w(p_{2})=q$, then by continuity of $w$ we could find two disjoint open
sets $U \subset \Omega(\mu)$, $p_{1} \in U$, and $V$, $p_{2} \in V$, such that
$w(U)=w(V)$. The density property of $\Omega(\mu)$ (see (4) in Definition \ref{nodedfamily})  asserts that there
are points $q_{1} \in U$ and $q_{2} \in V \cap \Omega(\mu)$ for which
$w(q_{1})=w(q_{2})$, what contradicts the one-to-one property of the map $w$ restricted to
$\Omega(\mu)$;
\item[(4)] $L^{\infty}(\Delta,G)_{1} \subset L^{\infty}_{{\rm noded}}(\Delta,G)$;
\item[(5)] As a consequence of (3) in Definition \ref{nodedbeltrami} we have $w(\widehat{\mathbb C})=\widehat{\mathbb C}$, that is, the map
$w$ is surjective.
 \end{itemize}
 \end{remark}

\s

The importance of the noded quasiconformal deformations of Kleinian groups is
reflected in the following result.

\s
\noindent
\begin{theorem}\label{Teo1} 
Let $\mu \in L^{\infty}_{{\rm noded}}(\Delta,G)$ be a 
noded Beltrami coefficient for $G$ supported in $\Delta$. 
If $w:\widehat{\mathbb C} \to \widehat{\mathbb C}$ is a 
noded quasiconformal deformation for $G$ with the complex dilatation $\mu$, then
there is a unique Kleinian group $\theta(G)$ of M\"obius transformations
and  a unique isomorphism of groups $\theta:G \to \theta(G)$ such
that $w g = \theta(g) w$.  Moreover, the region of discontinuity for the action of 
$\theta(G)$ on the Riemann sphere is $w(\Omega(\mu) \cap \Omega(G))$.
\end{theorem}

\begin{proof} Let $\mu$ be as in the hypothesis. For each $g \in G$ we 
proceed to obtain an orientation preserving homeomorphism $\theta(g)$ 
satisfying $w g = \theta(g) w$. Moreover by the construction, these
transformations are conformal automorphisms of the open set 
$w(\Omega(\mu))$. Let us denote the noded family of arcs associated with $\mu$
by ${\mathcal F}_{\mu}=\{\alpha_{1},....\}$ and let us set
$U=w(\widehat{\mathbb C}-\cup_{n}\alpha_{n})$.
Define $\theta(g):U \to U$ by the rule $w g = \theta(g) w$. It is well
defined because of the invariance property of $\Omega(\mu)$ given by
Proposition \ref{Prop1}  and the fact that $w$ restricted to 
$\widehat{\mathbb C}-\cup_{n}\alpha_{n}$
is a homeomorphism. 
The function $\theta(g):U \to U$ is a homeomorphism. 
We have that the set $w(\Omega(\mu))$
is invariant under the action of $\theta(g)$ for all $g \in G$. 
A direct computation asserts that $\theta(g)$ is a conformal automorphism of
$w(\Omega(\mu))$. For each arc $\alpha_{t} \in {\mathcal F}_{\mu}$ we have that
$w(\alpha_{t})$ is just a point $p_{t}$. The invariance property of
the arcs under the action of $G$ asserts that $g(\alpha_{n})$ is again
an arc $\alpha_{m}$ of the noded family:
we set $\theta(g)(p_{n})=p_{m}$. 

The following facts: (i) $w$ is continuous, (ii) $w$ restricted to $\Omega(\mu)$ is 
injective, (iii) $w$ preserves orientation,
(iv) $w(\Omega(\mu)) \cap w(\Lambda(\mu)) =\emptyset$, (v) the points
$p_{t}$ are all different and (vi) $w$ restricted to the complement of 
the arcs $\alpha_{t}$ is one-to-one, 
assert that the function $\theta(g)$ is an orientation preserving 
homeomorphism  of $\widehat{\mathbb C}$, such that $w g = \theta(g) w$.
Moreover, we have obtained a group $\theta(G)$, generated by the homeomorphisms
$\theta(g)$, and an isomorphism $\theta:G \to \theta(G)$ satisfying
$w g = \theta(g) w$.

The definition of the noded Beltrami coefficients for $G$ states 
the existence of a sequence $\{w_{n}\}$
of $\mu_{n}-$ quasiconformal homeomorphisms converging locally uniformly to $w$ in 
$\Omega(\mu)$ with the Beltrami coefficients $\mu_{n}$ for $G$. It follows 
that $w_{n} g w_{n}^{-1}$ is a M\"obius transformation for all $g \in G$. 
In particular, $\theta(g)$ is the local uniform limit of M\"obius 
transformations and, as a consequence, itself is a M\"obius transformation.
It follows that $\theta(G)$ is a group of M\"obius transformations with the
region of discontinuity $\Omega(\theta(G))=w(\Omega(\mu) \cap \Omega(G))$ as
claimed. \hspace{\fill}\end{proof}

\s

The arguments done at the end of the above proof give the following trivial result.

\s
\noindent
\begin{corollary}\label{Cor1}
Let $\mu \in \overline{L^{\infty}(\Delta,G)_{1}}$ 
and let $w:\Omega(\mu) \to w(\Omega(\mu))$ be a
homeomorphism with complex dilatation $\mu$.
Suppose that there is a sequence $w_{n}$ of quasiconformal homeomorphisms
with corresponding complex dilatations $\mu_{n} \in L^{\infty}(\Delta,G)_{1}$,
converging locally uniformly to $w$ in $\Omega(\mu)$. Then
there exist a group $\theta(G)$ of M\"obius transformations and 
an isomorphism of groups $\theta:G \to \theta(G)$, 
such that $w g = \theta(g) w$.
\end{corollary}

\s
\noindent
\begin{example}\label{Ex1}
Let $G$ be the cyclic group generated by the hyperbolic transformation
$C(z)=2z$  and let $\theta(G)$ be the cyclic group generated by the parabolic
transformation  $A(z)=z+1$. Define $w:\widehat{\mathbb C} \to \widehat{\mathbb C}$ as follows
$$ w(z)= \left\{\begin{array}{ll} h(|z|) + i g(\eta), & \mbox{ for } z=|z|e^{i\eta} \in 
         {\mathbb C}-[0,\infty],\\
          \infty,              & \mbox{ for } z \in [0,\infty], \end{array} 
\right. $$
where $h(x)=\ln(x)/\ln(2) -1/2$,\quad and \quad $g(\eta)=1/\tan(\eta/2)$. 
In this case ($z=|z|e^{i\eta}$) 
$$\mu(z)=(\frac{\mbox{$z$}}{\mbox{$\overline z$}})
 (\frac{\mbox{$2- \ln(2) \csc^2(\eta/2)$}}
      {\mbox{$2+\ln(2) \csc^2(\eta/2)$}}). $$

As we can see, the function $|\mu|$  depends only on the argument $\eta$.
The graphic of  $|\mu|$, as a function of the argument, is shown in 
Figure \ref{fig1}. We have $\|\mu\|_{\infty}=1$ and $\Lambda(\mu)=[0,+\infty]$. 
Moreover, it is not difficult to see that:
\begin{itemize}
\item[(a)] $\mu \in \overline{L^{\infty}(\Delta,G)_{1}}$;
\item[(b)] $w$ is a noded quasiconformal deformation of $G$; 
\item[(c)] $w g = \theta(g) w$, for every $g \in G$.
\end{itemize}

If we denote by $Q=\{1 \leq |z| \leq 2\}$ the fundamental annulus for $G$, then
a way to see the above is to consider for each $n \in \{1,2,3,...\}$:
\begin{itemize}
\item[(d)] the sequence of M\"obius transformations
$A_{n}(z)=\frac{n-1}{n+1}z+1$, which converges to
$A(z)$ as $n$ goes to $\infty$
\item[(e)] the sequence of the fundamental annuli $Q_{n}=\{ \frac{n-1}{2} \leq
|z-\frac{n+1}{2}|\leq \frac{n+1}{2} \}$ for the cyclic groups $G_{n}$ generated
by $A_{n}$; 
\item[(f)] a sequence of homeomorphisms $f_{n}:Q \to Q_{n}$, converging to
$w:Q \to R$, where $R=\{ 0 \leq \Re(z) \leq 1\}$, so that
$f_{n}(C(z))=A_{n}(f_{n}(z))$;
\item[(g)] the extension of $f_{n}$ by the relation
$f_{n}(C^{k}(z))=A^{k}_{n}(f_{n}(z))$;
\item[(h)] $\mu_{n}(z)=\frac{\overline{\partial}f_{n}(z)}{\partial
f_{n}(z)} \in L^{\infty}(\Delta,G)_{1}$, for $n=1,2,...$
\end{itemize}
In particular, $\mu \in L^{\infty}_{{\rm noded}}(\Delta,G)$. 

\begin{figure}
\centering
\includegraphics[width=15cm]{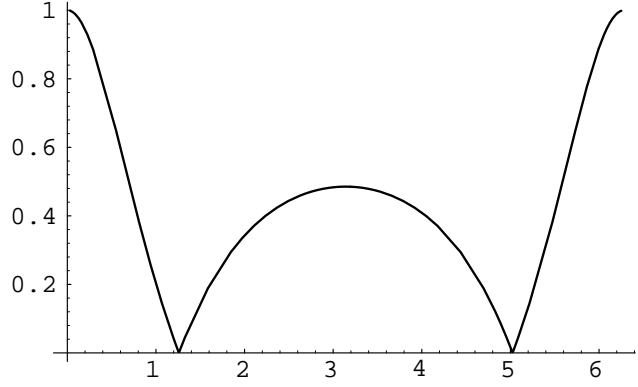}
\caption{Graph of $|\mu(e^{i\eta})|$}
\label{fig1}
\end{figure}
\end{example}

\s

For $\mu$ of Example \ref{Ex1} we have the following.  Assume 
that $w_{1}$ is another orientation preserving homeomorphism defined 
on $\Omega(\mu)$ solving the Beltrami equation for $\mu$. 
If we set $X=w_{1}(\Omega(\mu))$, then $T=w w_{1}^{-1}:X \to {\mathbb C}$
is a conformal homeomorphism. In particular, $X$ must be the complement 
of a point $p$ on the Riemann sphere. If we define $T(p)=\infty$ and $w_{1}|\gamma=p$, then the extension of $w_{1}$ is continuous and $T$ is necessarily
a M\"obius transformation. In particular, $w_{1}$ has a natural extension to a
noded quasiconformal deformation of $G$. For more general situation
for noded quasiconformal deformations we have the following lemma.

\s
\noindent
\begin{lemma}\label{Lema2}
Let $\mu \in \overline{L^{\infty}(\Delta,G)_{1}}$, for some
finitely generated Kleinian group $G$. Let $w_{1}$ and $w_{2}$ be noded
quasiconformal deformations of $G$ with the complex dilatation $\mu$. We can find
an orientation preserving homeomorphism  
$T:\widehat{\mathbb C} \to \widehat{\mathbb
C}$ with a conformal mapping $T:w_{1}(\Omega(\mu)) \to w_{2}(\Omega(\mu))$, such that
$T w_{1} = w_{2}$.
\end{lemma}

\begin{proof} The construction of $T$ is given as follows.
\begin{itemize}
\item[(3.1)] If $x \in w_{1}(\Omega(\mu))$, then set
$T(x)=w_{2}(w_{1}^{-1}(x))$.
\item[(3.2)] If $x \in w_{1}(\alpha_{n})$, then
set $T(x)=w_{2}(\alpha_{n})$.
\item[(3.3)] If $x \in
\Lambda(\mu)-\cup_{n}\alpha_{n}([0,1])$, then set $T(x)=w_{2}(w_{1}^{-1}(x))$.
\end{itemize}
\end{proof}

\s
\noindent
\begin{remark} 
Let $\mu \in L^{\infty}_{{\rm noded}}(\Delta,G)$ and let
$w:\widehat{\mathbb C} \to \widehat{\mathbb C}$ be  a
noded quasiconformal deformation of $G$ with complex dilatation $\mu$. If
$\theta(G)$ and $\theta:G \to \theta(G)$ are as in Theorem \ref{Teo1}, then
\begin{itemize} 
\item[(1)] if $p$ belongs to one of the arcs $\alpha_{i}$,
then  $w(p)$ is necessarily a doubly cusped parabolic fixed point;
\item[(2)] if $p \in \Lambda(G)$ is a loxodromic fixed point in $G$, 
which does not belong to any arc $\alpha_{i}$, then $w(p)$ is again a 
loxodromic fixed point;
\item[(3)] if $p$ is a rank two parabolic fixed point of $G$, then $w(p)$ is 
again a rank two parabolic fixed point in $\theta(G)$;
\item[(4)] if $p$ is a doubly cusped parabolic fixed point of $G$, then $w(p)$ 
is again doubly cusped in $\theta(G)$.
\end{itemize}
\end{remark}

\s
\noindent
\begin{theorem}\label{Teo2} 
Under the hypothesis of Theorem \ref{Teo1}, if $G$ is
geometrically finite and has an invariant component of its region of
discontinuity, then $\theta(G)$ is also geometrically finite.
\end{theorem}

\begin{proof} We use the same notations as in Theorem \ref{Teo1}. 
As we remarked above, we need only to check that if $x \in
\Lambda(\theta(G))$ is not a fixed point, then it is a point of
approximation. Let $p \in
\widehat{\mathbb C}$ be such that $x=w(p)$. 
Since $x$ is not a fixed point, the point $p$ can not  belong to any
arc $\alpha_{n}$. In the sequel, $w^{-1}(x)=\{p\}$ and $p \in \Lambda(G).$

Assume that $x$ is not a point of approximation in $\theta(G)$. Then for
every sequence $\{\theta(g_{n})\}$ of different elements of $\theta(G)$ there
is a point $z_{0} \in \widehat{\mathbb C}-\{x\}$ and a subsequence 
$\{\theta(g_{n_{k}})\}$, such that
$$d(\theta(g_{n_{k}})(x),\theta(g_{n_{k}})(z_{0})) \to 0,$$
where $d$ denotes the spherical diameter (see the definition of a point of
approximation given  in \cite[VI.B.1]{M1}). 

Let us choose a point $q \in \widehat{\mathbb C}-\{p\}$ with $w(q)=z_{0}$. 
The property $\theta(g) w = w g$ for all $g \in G$ yields that
$$d(w(g_{n_{k}}(p)),w(g_{n_{k}}(q))) \to 0.$$

We may assume, taking a suitable subsequence, that
$g_{n_{k}}(p)$ and $g_{n_{k}}(q)$ both converge, say to  $u$ and $v$,
respectively. Clearly, $u,v \in \Lambda(G)$. Then, the continuity of $w$ 
asserts that $$w(u)=w(v).$$

Since $p$ is a point of approximation for $G$, we may assume that for
the sequence $\{g_{n}\}$  there exists $\delta > 0$, such
that $$d(g_{n}(p),g_{n}(z)) \geq \delta > 0,$$
uniformly in compact sets of $\widehat{\mathbb C}-\{p\}$. In particular, we
must have that $$u \neq v.$$

It follows that there is $n$, such that $u,v \in \alpha_{n}$. The
facts that $u \neq v$ and $u,v \in \Lambda(G)$, imply that $u$ and $v$ are
two extremes of the same arc $\alpha_{n}$, i.e., the fixed points
of a loxodromic element $g \in G_{n}<G$, where $G_{n}$ is the stabilizer of
$\alpha_{n}$.

Since there is an invariant component $\Delta \subset \Omega(G)$,
we may construct a simple arc $\gamma$ starting at $p$ and ending at $q$,
so that $\gamma -\{p,q\} \subset \Delta^{0}$, where $\Delta^{0}$ denotes
$\Delta$ minus elliptic fixed points.

Consider two open balls $B_{1}$ and $B_{2}$ in $\Delta$
containing no elliptic fixed point, such that $B_{1}$ covers one extreme
of  $\gamma-\{p,q\}$ (e.g., the extreme determined by $p$) and $B_{2}$ covers
the other extreme. We can use compactness arguments to find a finite number 
of balls $B_{3}$,..., $B_{r}$, in $\Delta$ each, containing no
elliptic fixed point, and 
$\gamma-\{p,q\} \subset B_{1} \cup \cdots \cup B_{r}$.

Proposition B.3 in \cite[page 17]{M1} asserts that the spherical diameter of
$g_{n}(B_{i})$ goes to zero as $n$ goes to infinity. In particular,
the spherical diameter of $\gamma_{n}=g_{n}(\gamma)$ goes to zero as
$n$ goes to infinity. But the diameter of the sequence $g_{n_{k}}(\gamma)$
goes to a value bigger or equal than  
the spherical distance between $u$ and $v$ as $n$ goes to $\infty$, that 
contradicts the inequality $u \neq v$.
\end{proof}

\s

In Theorem \ref{Teo2} we have used the extra hypothesis that $G$
is a function group, that is, we have assumed the existence of an invariant component $\Delta$ in its region of discontinuity. This condition was needed at the end of
the proof in order to find an arc $\gamma \subset \Omega(G)$ connecting the
points $p$ and $q$ obtained in the proof. We may avoid such
an extra condition by the following argumentation. We may think of $G$ as a group of
conformal automorphisms of the $3$-dimensional
sphere $S^{3}$. We can now find an arc $\gamma \subset \Omega_{S^{3}}(G)$
connecting the limit points $p$ and $q$. Then we can use the fact that 
the spherical diameter of $g_{n}(B)$ goes to zero as $n$ goes to $\infty$
for each Euclidean ball $B \subset \Omega_{S^{3}}(G)$. 
This allows to rewrite Theorem \ref{Teo2} as follows.

\s
\noindent
\begin{theorem}\label{Teo3} 
Under the hypothesis of Theorem \ref{Teo1}, if $G$ is
geometrically finite, then $\theta(G)$ is also geometrically finite.
\end{theorem}

\section{Topological Realizations}
Let us consider a Kleinian group $G$ and a $G$-invariant
collection $\Delta$ of components of its region of discontinuity $\Omega(G)$.

A simple loop 
$\gamma \subset S$ (or a simple arc connecting two branch values of order $2$) is called pinchable if (a) its lifting on $\Omega(G)$
consists of pairwise disjoint simple arcs, and (b) each of these components is
stabilized by a primitive loxodromic transformation in $G$.

The stabilizer of a component of the lifting of a pinchable
loop $\gamma$ is a cyclic group generated by a primitive loxodromic
transformation, and the stabilizer of a component of the lifting of a
pinchable arc is a ${\mathbb Z}/2{\mathbb Z}$-extension of a cyclic group
generated by a primitive loxodromic transformation. We say that such a cyclic
group (or ${\mathbb Z}/2{\mathbb Z}$-extension) is defined by the pinchable
loop $\gamma$.

If $\alpha$ and $\beta$ are two components of the lifting of a pinchable
$\gamma$, then their stabilizers are conjugate in $G$. We say that a
collection $\{\gamma_{1},...,\gamma_{m}\}$ of pairwise disjoint pinchable
loops or arcs is admissible if  they define non-conjugate 
groups in $G$ for $i \neq j$.

For each collection ${\mathcal F}=\{\gamma_{1},...,\gamma_{m}\}$ of pinchable
loops on $S$ we define the following equivalence relation on the Riemann sphere
$\widehat{\mathbb C}$. Two points $p, q \in \widehat{\mathbb C}$
are equivalent if either:
\begin{itemize}
\item[(1)] $p=q$; or
\item[(2)] there is a component $\widetilde{\gamma_{j}}$ of the lifting of
some pinchable loop or arc $\gamma_{j}$, such that $p,q \in
\widetilde{\gamma_{j}} \cup \{a,b\}$, where  $a$ and $b$ are the endpoints of
$\widetilde{\gamma_{j}}$ (that is, the fixed points of a primitive
loxodromic transformation in the stabilizer of
$\widetilde{\gamma_{j}}$ in $G$)  
\end{itemize}

We claim that the set of equivalence classes for such an equivalence relation is
(topologically) the Riemann sphere. In fact, let us denote by $\widetilde{\mathcal F}$ the collection of all arcs (included their endpoints), as considered above in (2), and let us consider the collection of continua given by the collection of arcs in $\widetilde{\mathcal F}$ as points and also each of the points in the complement of $\widetilde{\mathcal F}$. The discreteness of $G$ asserts that such collections of points are  semi-continuous as defined in \cite{Moore}. Now the result  follows from Theorem 22 from \cite{Moore}.

Let us denote by
$P:\widehat{\mathbb C} \to \widehat{\mathbb C}$ the natural continuous
projection defined by the above relation. 

\s
\noindent
\begin{question} Is there some orientation preserving homeomorphism
$Q:\widehat{\mathbb C} \to \widehat{\mathbb C}$, such that the map
$Q P:\widehat{\mathbb C} \to \widehat{\mathbb C}$ is a noded quasiconformal
deformation of the group $G$? 
\end{question}

\s

The same arguments used in the proof of Theorem \ref{Teo1} permit us to obtain the
following.

\s
\noindent
\begin{proposition} \label{Prop3}
Let $G$ be a Kleinian group and ${\mathcal F}$ be a
collection of admissible pinchable loops. Denote by 
$P:\widehat{\mathbb C} \to \widehat{\mathbb C}$ the continuous projection naturally  induced by the equivalence relation defined by $\mathcal F$. Then we
have a Kleinian group $\theta(G)$ of orientation preserving homeomorphisms of
the Riemann sphere and an isomorphism $\theta:G \to \theta(G)$, such that $P g =
\theta(g) P$, for all $g \in G$.
\end{proposition}

\s
\noindent
\begin{remark} 
Theorem \ref{Teo2} (and also the results of \cite{M4} and \cite{Ohshika})
asserts that if $G$ is geometrically finite, then 
$\theta(G)$ is also geometrically finite.
\end{remark}

\section{Construction of Noded Quasiconformal Deformations} \label{ejemplo}
In this section we are aimed at construction of a sequence of quasiconformal
maps  $f_n:\widehat{\mathbb C}\to \widehat{\mathbb C}$ that converges to a
noded quasicoformal deformation $w$ locally uniformly almost everywhere
in $\widehat{\mathbb C} \setminus  \Lambda({\mu})$. Each $f_n$ satisfies the
Beltrami  equation
$\bar{\partial}f_n(z)=\mu_n(z)\partial f_n(z)$ with 
$\mu_n\in L^{\infty}(\Delta,G)_1$.
To avoid actions of the M\"obius group in $\widehat{\mathbb C}$ we require  
certain normalization
for all maps $f_n$ assuming 
$f_n(z)=z+a_0+a_1/z+\dots$ for $z$ close to $\infty$. The same normalization we assume for $w$. 
We can suppose
that the maps $f_n$ and $w$ are conformal in $|z|>R$ for some $R$ sufficiently 
big, such that
$$\mathcal F_{\mu}=\{\alpha_{1},...\}\subset \{|z|<R\}.$$ 

We also define a
sequence of the domain systems in the following way. For an arc $\alpha_k$
there exists  a simply connected domain $\Delta_{k,1}$ and a conformal map
$Z=g_{k}(z)$ with the following properties. The domain $\Delta_{k,1}$ contains
the curve $\alpha_k$ so that the endpoints $a_k,b_k$ of $\alpha_k$ lie in its
border, $\Delta_{k,1}\bigcap\Delta_{j,1}=\emptyset$, $j\neq k$, moreover,
$G(\Delta_{k,1})$ is a system of non-overlapping domains. We mark four
points 1,2,3,4 at the border of $\Delta_{k,1}$ so that $a_k$ lies on the side
1,2 of the obtained curved quadrangle and $b_k$ lies on the side 3,4. The
rectifying conformal map $Z=g_k(z)$ transforms  the  quadrangle $\Delta_{k,1}$
into the rectangle $R_{k,1}$ in $(Z)$-plane with the vertices
$-i/2,i/2,l_k+i/2,l_k-i/2$, $g_k(a_k)=0$, $g_k(b_k)=l_k$,
$g_k(\alpha_k)=[0,l_k]$ with the obvious correspondence of the vertices,
$1/l_k$ is the modulus of the quadrangle $\Delta_{k,1}$ with respect to the 
family of curves that connect the sides 1,2 and 3,4. Now we set a
quasiconformal map
$W=\omega_{k,1}(Z)\equiv\frac{Z-\delta_1\bar{Z}}{1+\delta_1}$, $\delta_1\in
[0,1)$ that maps $R_{k,1}$ onto the rectangle $R'_{k,1}$ in the $(W)$-plane
preserving three points $\omega_{k,1}(0)=0,\omega_{k,1}(\infty)=\infty,
\omega_{k,1}(i)=i$. So this is a unique $\delta_1$-quasiconformal map
$\omega_{k,1}:\widehat{\mathbb C} \to \widehat{\mathbb C}$ having the
complex dilatation $(-\delta_1)$ with such a normalization.

Choose the domains
$\Delta_{k,n}=g^{-1}_k(R_n)$, where $R_{k,n}$ is the rectangle in $(Z)$-plane
with the vertices $-i/2n,i/2n,l_k+i/2n,l_k-i/2n$

Now we define a sequence of the Beltrami
coefficients $\mu_n(z)$, $z\in \widehat{\mathbb C}$ in the following way. 
Let $\mu(z)$, $z\in \widehat{\mathbb C}$ be a noded
Beltrami coefficient  degenerating ($|\mu(z)|=1$) in $\mathcal F_{\mu}$, 
and $\mu(z)=0$ in
$|z|>R$ (by the normalization of the functions $f_n$ and $w$).
Assume that it is a continuous function in $\widehat{\mathbb C}$ and 
$|\mu(z)|_{z\in \partial\,\Delta_{k,n}}=\delta_n=1-1/(n+1)^4$
(independently of $k$, $n=1,2,\dots$). We put $\mu_n(z)=\mu(z)$ everywhere in
$\widehat{\mathbb C}-\bigcup_{k=1}^{\infty}\Delta_{k,n}$. Then we
assume $\mu_n(z)=-\delta_{n}e^{i\theta_{n,k}(z)}$,
$\theta_{n,k}(z)=\arg\,(\overline{g'_k(z)}/g'_k(z))$, $z\in \Delta_{k,n}$ for the chosen
$k$, and  extend $\theta_{n,k}(z)$ for other $k=1,2,\dots$ by
the actions of the Kleinian group $G$. Construct a sequence $w=f_n(z)$ of
normalized quasiconformal maps satisfying the Beltrami equation with the
complex dilatation $\mu_n$. By the construction of $\mu_n$, there is a
conformal map $W=h_{k,n}(w)$ of the quadrangle $f_n(\Delta_{k,n})$ onto the
rectangle $R_{k,n}'=\omega_{k,n}(R_{k,n})$, such that $$h_{k,n}^{-1}\circ
\omega_n\circ g_k\bigg|_{\Delta_{k,n}}\equiv f_n\bigg|_{\Delta_{k,n}}.$$ 

Here
one can obtain the map $\omega_{k,n}$ substituting $\delta_1$ by $\delta_n$ in
$\omega_{k,1}$. Denote by $|\cdot|$ the Euclidean length. 

We have 
$$|f_n(\alpha_k)|=\int\limits_{f(\alpha_k)}|dw|=$$
$$\int\limits_{\alpha_k}|\partial\,f_n(z)||1+\mu_n(z)d\bar{z}/dz||dz|
=\int\limits_{0}^{l_k}\bigg|\frac{\partial}{\partial\, Z}(f_n\circ
g^{-1}_k)\bigg|(1-\delta_n)|dZ|. $$

Integrating along the imaginary axis in the $(Z)$-plane and using the Fubini
formula we deduce that
$$
\frac{1}{n}|f_n(\alpha_k)|=\iint\limits_{R_{k,n}}\bigg|
\frac{\partial}{\partial\, Z}(f_n\circ g_k^{-1})\bigg|(1-\delta_n)d\sigma_Z,
$$ where $d\sigma_Z$ denotes the area element.

Then using the Schwarz inequality we derive 
$$
\frac{1}{n}|f_n(\alpha_k)|\leq\left(\iint\limits_{R_{k,n}}
\bigg|\frac{\partial}{\partial\, Z}(f_n\circ
g^{-1}_k)\bigg|^{2}(1-\delta_n^2)d\sigma_Z\right)^{\frac{1}{2}}
\sqrt{\frac{1-\delta_n}{1+\delta_n}}\sqrt{\frac{l_k}{n}}=$$ $$
=\sqrt{S_{k,n}}\left(\frac{l_k}{n(1+(1+n)^4)}\right)^{\frac{1}{2}}, 
$$
where $S_{k,n}$ is the Euclidean area of $f_n(\Delta_{k,n})$. 
The area distortion theorem 
for univalent functions with the normalization $z+a_0+a_1/z+\dots$, for
$|z|>R$, implies the estimate  $S_{k,n}<\pi (1+R)^4$. So, 
$$
|f_n(\alpha_k)|<(1+R)^{2}\left(\frac{\pi n
l_k}{1+(1+n)^4}\right)^{\frac{1}{2}}\to 0, 
\,\,\text{as}\,\,n\to\infty.
$$
Thus, we construct the sequence of quasiconformal maps that shrinks 
the curves
from $\mathcal F_{\mu}$ into points.

\section{The Noded Teichm\"uller space of $G$ supported in $\Delta$}\label{construccion}
Let $G$ be a finitely generated Kleinian group and $\Delta$ be a collection of components of its region of discontinuity, invariant under the action of $G$.
We are in conditions to define a partial closure of $T(G,\Delta)$ in a natural sense.
As a consequence of Theorems \ref{Teo2} and \ref{Teo3}  and Lemma \ref{Lema2}, we may extend the Teichm\"uller equivalence relation, given in the first section on  $L^{\infty}(\Delta,G)_{1}$, to the whole $L^{\infty}_{{\rm noded}}(\Delta,G)$ as follows. Let $\mu$ and $\nu$ be in $L^{\infty}_{{\rm noded}}(\Delta,G)$ and $w_{\mu}$, $w_{\nu}$ be associated noded quasiconformal deformations for $G$, respectively. 
Theorems \ref{Teo2} and \ref{Teo3} assert the existence of
isomorphisms $$\theta_{\mu}:G \to G_{\mu} \quad \mbox{and} \quad 
\theta_{\nu}:G \to G_{\nu}\;,$$ where 
$G_{\mu}$ and $G_{\nu}$ are Kleinian groups 
such that $w_{\mu} g = \theta_{\mu}(g) w_{\mu}$ and 
$w_{\nu} g = \theta_{\nu}(g) w_{\nu}$  for all $g \in G$. 
We say that $\mu$ and $\nu$ are {\it noded Teichm\"uller
equivalent} if there is an orientation preserving homeomorphism 
$A:\widehat{\mathbb C} \to \widehat{\mathbb C}$, such that 
\begin{itemize}
\item[(1)] $A(w_{\mu}(\Omega(\mu)))=w_{\nu}(\Omega(\nu))$;
\item[(2)] $A:w_{\mu}(\Omega(\mu)) \to w_{\nu}(\Omega(\nu))$ is conformal;
\item[(3)] $\theta_{\nu}(g)=A \theta_{\nu}(g) A^{-1}$ for all $g \in G$.
\end{itemize}

We call the set of equivalence classes of noded Beltrami coefficients
for $G$ supported in $\Delta$, denoted by $NT(\Delta,G)$,  the
{\it noded  deformation space of $G$ supported in $\Delta$}. If $\Delta=\Omega(G)$,
then we denote it by $NT(G)$.

\s
\noindent
\begin{remark} If $\mu, \nu \in L^{\infty}_{{\rm noded}}(\Delta,G)$ are noded
Teichm\"uller equivalent and $\mu \in L^{\infty}(\Delta,G)_{1}$, then (1) and
(2) in above definition assert that $\nu \in L^{\infty}(\Delta,G)_{1}$ and
that they are Teichm\"uller equivalent. Moreover,
the inclusion $L^{\infty}(\Delta,G)_{1} \subset L^{\infty}_{{\rm noded}}(\Delta,G)$ induces, under the above equivalence relation, 
the inclusion $$T(\Delta,G) \subset NT(\Delta,G).$$
\end{remark}

\s

In the next we consider some particular classes of finitely generated Kleinian groups to describe $NT(G,\Delta)$ and relate it to already known partial closures of Teichm\"uller space.

\s
\noindent
\begin{example}\label{Ex2}
If $G$ is either a finite Kleinian group or a
finite extension of a purely parabolic Kleinian group, then we have 
$$L^{\infty}_{{\rm noded}}(\Omega(G),G)=L^{\infty}(\Omega(G),G)_{1}=
L^{\infty}(\widehat{\mathbb C})_{1},$$
in particular, $NT(G)=T(G)$. This example shows the importance of loxodromic elements in $G$ in order two get a non-trivial partial closure.
\end{example}

\s
\noindent
\begin{example}\label{Ex3}
Let $G$ be the cyclic loxodromic group generated by $C(z)=2z$. In this case,
we necessarily have that $\Delta=\Omega(G)={\mathbb C}-\{0\}={\mathbb C}^{*}$.
Example \ref{Ex1} asserts that 
$$L^{\infty}_{{\rm noded}}({\mathbb C}^{*},G)-L^{\infty}({\mathbb C}^{*},G)_{1} 
\neq \emptyset.$$

Let $\mu$ and $\nu$ be in 
$L^{\infty}_{{\rm noded}}({\mathbb C}^{*},G)-L^{\infty}({\mathbb C}^{*},G)_{1}$. 
Each of the associated noded families of arcs consists exactly of one arc. 
Denote these noded families by
${\mathcal F}_{\mu}=\{\alpha_{\mu}\}$ and ${\mathcal F}_{\nu}=\{\alpha_{\nu}\}$.
Assume that $w_{\mu}:\widehat{\mathbb C} \to \widehat{\mathbb C}$ and 
$w_{\nu}:\widehat{\mathbb C} \to \widehat{\mathbb C}$ are corresponding
noded quasiconformal deformations for $G$ with the complex dilatations $\mu$ and 
$\nu$, respectively.
Set $w_{\mu}(\alpha_{\mu})=p_{\mu}$ and 
$w_{\nu}(\alpha_{\nu})=p_{\nu}$. 
We have the natural isomorphisms $\theta_{\mu}:G \to G_{\mu}$ and 
$\theta_{\nu}:G \to G_{\nu}$, where $G_{\mu}$ and $G_{\nu}$ are cyclic
groups of orientation preserving homeomorphisms of $\widehat{\mathbb C}$, 
such that
$w_{\mu} C = \theta_{\mu}(C) w_{\mu}$ and 
$w_{\nu} C = \theta_{\nu}(C) w_{\nu}$.
We have that $\theta_{\mu}(C)$ and $\theta_{\nu}(C)$ are conformal 
homeomorphisms of $\widehat{\mathbb C}-\{p_{\mu}\}$ and 
$\widehat{\mathbb C}-\{p_{\nu}\}$,
respectively. In particular, they are M\"obius transformations. Since they
have exactly one fixed point, they are parabolic and analytically 
conjugate. The conjugacy is given by the M\"obius transformation 
$T:\widehat{\mathbb C} \to \widehat{\mathbb C}$ defined by 
$T(p_{\mu})=p_{\nu}$, and  $T=w_{\mu}^{-1}w_{\nu}$ on the complement.  
This transformation satisfies the equality 
$\theta_{\nu}(C) = T \theta_{\mu}(C) T^{-1}$. 
Now we have that $\mu$ and $\nu$ are noded Teichm\"uller equivalent. It follows
that the partial closure $NT(G)$ of the deformation space $T(G)$ (the Schottky
space of genus one, which is holomorphically equivalent to the punctured unit
disc) is given by adding a point to the boundary of it. This added point corresponds exactly to the point we need to add to the genus one moduli space to get Deligne-Mumford's compactification; which is the Riemann sphere.
\end{example}

\s
\noindent
\begin{example}\label{Ex4}
Let $G$ be a Schottky group of genus $g \geq 2$.
Set $\Omega$ its region of discontinuity. 
In this case, we have that $L^{\infty}_{{\rm noded}}(\Omega,G)$ contains more 
points than $L^{\infty}(\Omega,G)_{1}$. This is just a consequence of the
construction done in Section 4 and also from the fact that Example \ref{Ex3} can be
easily generalized to any genus. For each $\mu \in L^{\infty}_{{\rm noded}}(\Omega,G)$, the Kleinian group $\theta(G)$ of Theorem \ref{Teo1} is isomorphic to $G$ and, by Theorem \ref{Teo2}, geometrically finite. It follows from the results of \cite{H1} that $\theta(G)$ is a noded Schottky group of rank $g$. It is not difficult to see that each noded Schottky group of rank $g$ may be obtained in the above way. In particular, we have that $NT(G)$ is exactly the noded Schottky space of \cite{H1}. In \cite{H-M} we have observed that most of the noded Schottky groups can not be defined by circles, in particular, the noded Schottky groups (which are not Schottky groups) obtained by allowance of tangencies of Schottky circles as in \cite{G-W} only form a small portion of the boundary of noded Schottky space. 
\end{example}

\s
\noindent
\begin{example}\label{Ex5}
Let $G$ be a torsion free co-compact Fuchsian
group, keeping the unit disc $\Delta_{1}$ invariant. 
Set $\Delta_{2}=\widehat{\mathbb C}-\overline{\Delta_{1}}$.
We have that $T(\Delta_{i},G)$ is the Teichm\"uller space of $G$, which is a 
simply connected complex manifold of complex dimension $3g-3$ (here $g$ is 
the genus of $G$). Another space is $T(G)$, the deformation space of $G$. This
is a simply connected complex manifold of complex dimension $6g-6$. 
We have that $T(G)$ is complex holomorphically equivalent to
$T(\Delta_{1},G) \times T(\Delta_{2},G)$. This equivalence is given
by the holomorphic homeomorphism
$$[\mu] \to ([\mu_{1}],[\mu_{2}]),$$
where $\mu_{i}$ is defined as $\mu$ on $\Delta_{i}$ and zero on its
complement. We also have the partial closures 
$NT(G)$, $NT(\Delta_{1},G)$, and $NT(\Delta_{2},G)$. The above holomorphic map
can be extended continuously to a homeomorphism between $NT(G)$ and 
$NT(\Delta_{1},G) \times NT(\Delta_{2},G)$.
Each $NT(\Delta_{i},G)$ can be identified with Abikoff's Augmented Teichm\"uller
space of $G$ defined in \cite{Ab1} and \cite{Ab2}. Related to this case, see
also \cite{H2} and \cite{Kra-Maskit}.
\end{example}

\s
\noindent
\begin{remark} In the setting of representation of groups, 
we consider in the representation space
Hom($G$,PGL($2,{\mathbb C}$)) those faithful representations of $G$ which are 
geometrically represented by noded quasiconformal deformations. 
This is a natural generalization for the deformation space of $G$,
on which one considers the geometric representations given by
quasiconformal homeomorphisms of $G$.
\end{remark}


\end{document}